\documentclass[square]{ws-procs975x65}

\usepackage{wrapfig}
\input labelfig.tex

\def\Z{{\mathbb Z}}
\def\R{{\mathbb R}}

\def\Quad{\mathcal{Q}}
\def\eps{\epsilon}

\def\Grimm{Grimmett:newbook}
\def\HPSIIC{HPS:IIC}
\def\WWperc{arXiv:0710.0856}
\def\ABNW{MR2001c:60151}
\def\SmirnovPerc{MR1851632}
\def\DPSL{DPSL}
\def\SchrammSmirnovNoise{SSblacknoise}
\def\KestenScaling{MR88k:60174}
\def\GPS1{arXiv:0803.3750}
\def\SteifSurvey{arXiv:0901.4760}
\def\OdedSLE{MR2001m:60227}
\def\RGB{Wilson:RGB}
\def\CFN{MR2284886}
\def\LPSMSF{MR2271476}
\def\DamronSV{arXiv:0806.2425}
\def\Alexander{MR96c:60114}

\begin{document}

\title{The scaling limit of the Minimal Spanning Tree --- a preliminary report}

\author{G\'abor Pete}

\address{Department of Mathematics, University of Toronto, Canada\\
http://www.math.toronto.edu/\~{}gabor}

\author{joint work with Christophe Garban$^*$ and Oded Schramm$^{**}$}
\address{$^*$ENS Lyon\\
$^{**}$Microsoft Research (died on Sept.~1, 2008)}

\begin{abstract}
This is a brief description of how the recent proof of the existence and conformal covariance of the scaling limits of dynamical and near-critical planar percolation implies the existence and several topological properties of the scaling limit of the Minimal Spanning Tree, and that it is invariant under scalings, rotations and translations. However, we do not expect conformal invariance: we explain why not and what is missing for a proof.
\end{abstract}


\bodymatter

\section{Introduction}\label{s.intro}

Critical planar percolation, both the discrete process and the continuum scaling limit, have become central objects of probability theory and statistical mechanics; see \cite{\Grimm} for the classical results and \cite{\WWperc} for a great course on our present knowledge, using conformal invariance \cite{\SmirnovPerc} and SLE \cite{\OdedSLE}.
In the past few years, there has also been a lot of progress on dynamical percolation, which is not only the natural time evolution with critical percolation as the stationary measure, and the natural framework to study how noise effects the system and how it produces exceptional events, but also provides tools to understand the near-critical regime and related objects like Invasion Percolation and the Minimal Spanning Tree (MST). See the survey \cite{\SteifSurvey}, which discusses not only the recent \cite{\GPS1}, but also some work in preparation \cite{\DPSL,\HPSIIC}. In this note, we explain the applications to the MST, partly to encourage others to work on the main remaining open problem: conformal non-invariance of the MST scaling limit, which is certainly interesting given the translational, rotational and scale invariance that we can now prove.

We thank David Wilson for conversations and one of the pictures. The work of GP was supported by an NSERC Discovery Grant.

\section{Dynamical and near-critical percolation}\label{s.DPNCE}

In a series of papers \cite{DPSL}, we show that, in dynamical percolation, if each site of the triangular grid with mesh $\eta$ has a Poisson clock with rate $r(\eta)=\eta^2\alpha_4(\eta,1)=\eta^{3/4+o(1)}$ switching between black and white, then the $\eta\downarrow 0$ limit of this system exists as a Markov process. If the clocks always switch from white to black when they ring, then started from critical percolation at time 0, at time $t$ we have near-critical percolation with density roughly $1/2+t r(\eta)$ for black, for $t\in(-\infty,\infty)$. We show that this near-critical ensemble also has a scaling limit (called NCESL). We should briefly mention here what the topological space is where these limit processes live: following the description of \cite{\SchrammSmirnovNoise} for the scaling limit of static critical percolation, this is a compact metrizable space that encodes all macroscopic crossing events of ``quads'' (conformal rectangles), i.e., an element of the space tells which quads are crossed by the percolation configuration. Then, partly following the suggestion of \cite{\CFN}, we build the two limits from critical percolation, in two main steps: 

{\bf (1)} The {\bf normalized counting measure} on the sites that are pivotal for the left-right crossing of any given quad $\Quad$ converge to a finite measure $\mu^\Quad$ that is measurable in the scaling limit of critical percolation. We also show that this measure is {\bf conformally covariant}: if the domain is changed by $\phi(z)$, then we get $\mu^{\phi(\Quad)}$ from $\mu^\Quad$ by scaling locally by $|\phi'(z)|^{3/4}$: there are more pivotals for a larger domain. Finally, we show that the collection of these pivotal measures $\mu^\Quad$ can be used to understand the ``importance measures'' $\mu^\eps$: the amount of sites that have the alternating 4-arm event to macroscopic distance at least $\eps$.

{\bf  (2)} {\bf Stability:} Fix a quad $\Quad$, and let the set of sites switched in $[0,t]$ be $W_t$. Then the probability that a configuration $\omega_\eta$ can be changed on $W_t$ into $\omega_\eta',\omega_\eta''$ such that they agree on any site that is at least $\eps$-important in $\omega_\eta$, but $\Quad$ is crossed by $\omega_\eta'$ while not crossed by $\omega_\eta''$, is small if $\eps$ is small, uniformly in the mesh $\eta$. Note that this is a strengthening of Kesten's theorem \cite{\KestenScaling} that the 4-arm probabilities remain comparable in the entire near-critical regime. 

By (1), using only macroscopic information, we can tell for all $\eps>0$ how many microscopic $\eps$-important sites there are in any region, and hence we know the rate with which important switches start happening when we start the dynamics. Then part (2) says that by following the switches of all these initially $\eps$-important sites we can predict well (as $\eps\downarrow 0$) the state of any quad crossing event at any later time. Hence we get a well-defined Markov process in the scaling limit, both in the dynamical and the near-critical cases. We also get that these scaling limits are {\bf conformally covariant}: time is scaled locally by $|\phi'(z)|^{3/4}$.

These results have (or may have) several applications. Similarly to the measurable measure on pivotals, we can construct the limit of the length measure on a percolation interface, giving the first {\bf physical  time-parametrization} for the SLE$_6$ curve (as opposed to the conformal capacity parametrization). We hope to describe {\bf near-critical interfaces} with a {\bf massive $SLE_6$}, involving a self-interacting drift term. The rotational invariance of the NCESL seems to help prove the asymptotic circularity of the percolation {\bf Wulff crystal}, as $p\downarrow p_c$. We are also planning to study the dynamical and near-critical {\bf FK-Ising models}. However, there are limitations to our methods: we do not have any {\bf near-critical Cardy's formula}, and do not have a guess for the {\bf dimension} of Minimal Spanning Tree paths.

\section{The Minimal Spanning Tree}\label{s.MST}

\begin{wrapfigure}[8]{r}{1.5 in}
\begin{center}
\vspace{-0.5 in}
\epsfysize=1.5 in \epsffile{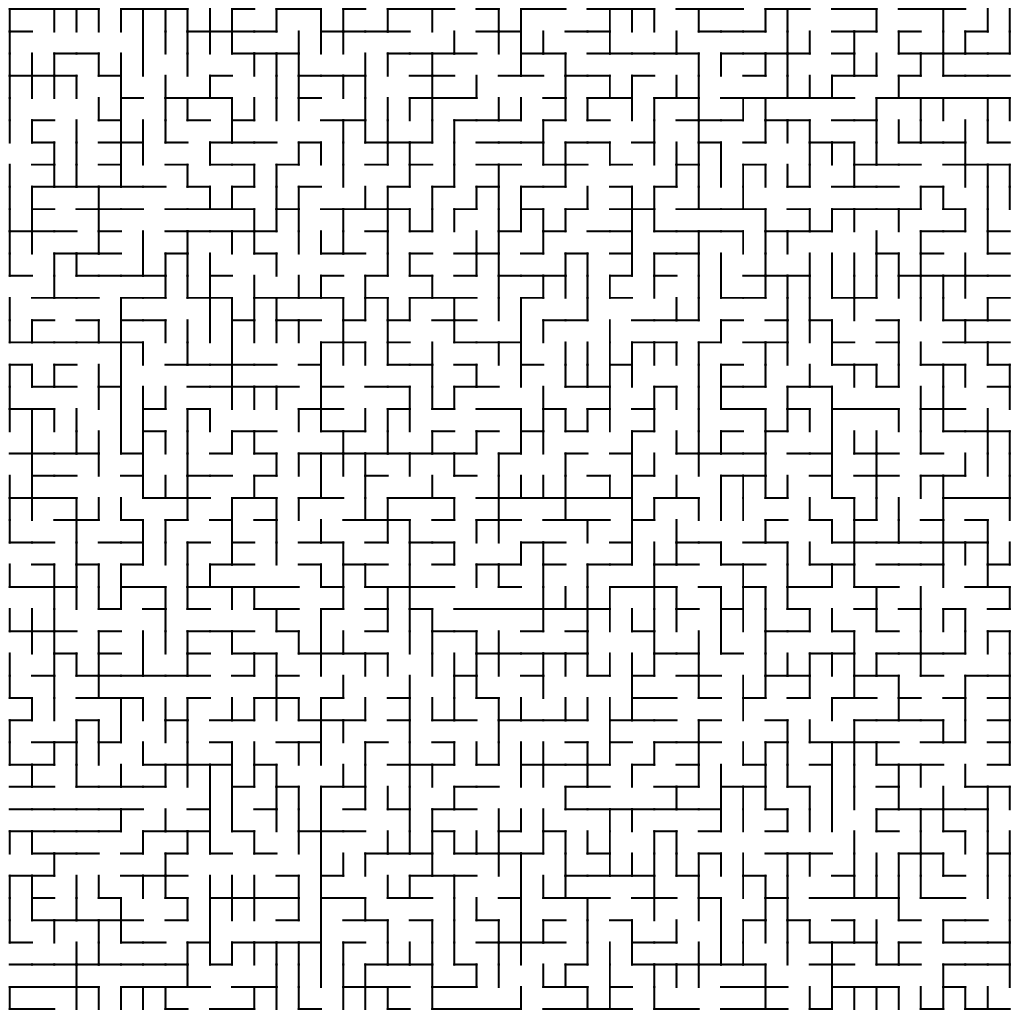}
\end{center}
\end{wrapfigure}

For each edge of a finite graph, say $e\in E(\Z_n^2)$, let $U(e)$ be
an i.i.d.~Unif$[0,1]$ label. The {\bf Minimal Spanning Tree} is the spanning tree $T$
for which $\sum_{e\in T} U(e)$ is minimal. This is well-known to be the same
as the collection of lowest level paths between all pairs of vertices (i.e., 
the path between the two points for which the maximum label on the path is minimal).
Or, delete from each cycle the edge with the highest label $U$. 
This also shows that $T$ depends only on the ordering of the labels, not the values themselves.

We can also use Unif$[0,1]$ labels to get a coupling of percolation for all densities $p$, and use the same labels to get the MST. This way we get a {\bf coupling} between the MST and the {\bf percolation ensemble}. Moreover, the macroscopic structure of the MST is basically determined by the labels in the near-critical regime, as follows.

\begin{wrapfigure}[6]{r}{1.4 in}
\begin{center}
\vspace{-0.35 in}
\epsfxsize=1.4 in \epsffile{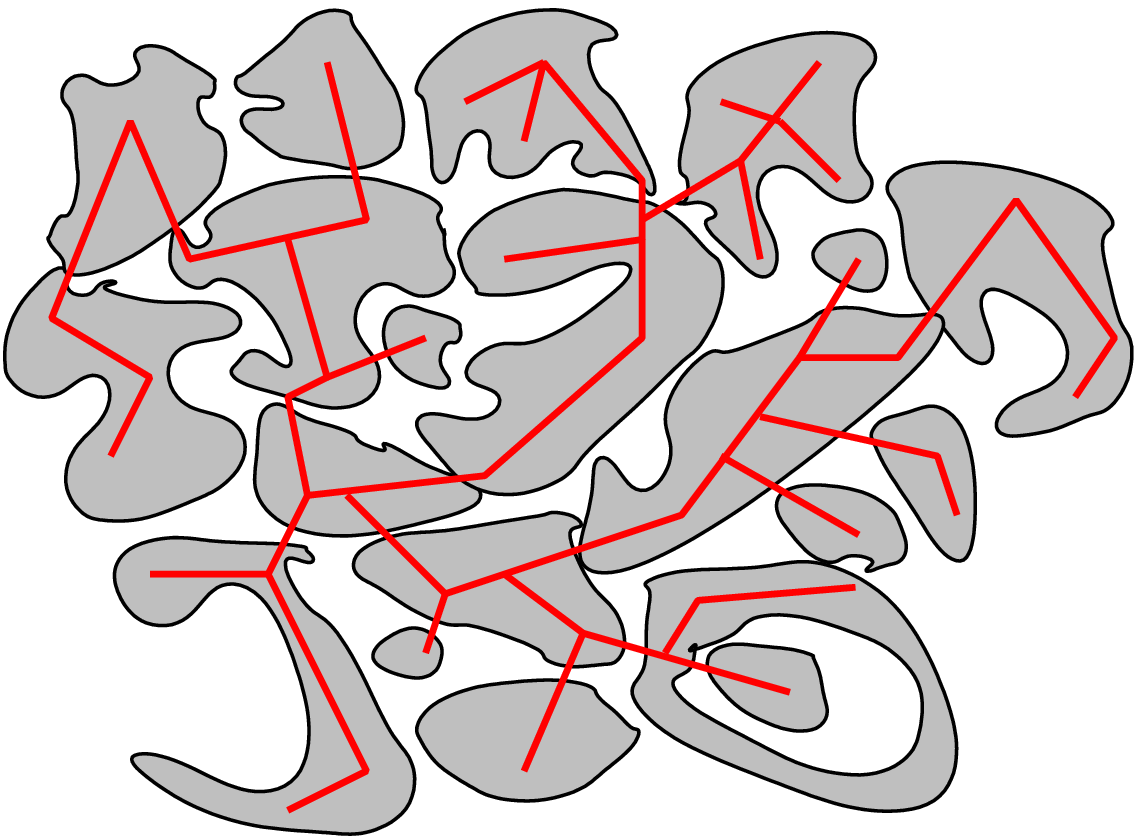}
\end{center}
\end{wrapfigure}

Consider the {\bf $\lambda$-clusters} in NCE, for $\lambda\in(-\infty,\infty)$, i.e., the connected components given by the labels at most $1/2+\lambda\, r(\eta)$. Contract each component into a single vertex, resulting in the ``cluster graph''. It is easy to verify that making these contractions on the MST we get exactly the MST on the cluster graph. We denote this {\bf cluster tree} by $T_\lambda$. Since the largest $\lambda$-clusters for $\lambda \ll 0$ are of very small macroscopic size, the tree $T_\lambda$ will tell us the macroscopic structure of the MST. On the other hand, for $\lambda\gg 0$, most sites are in a few large $\lambda$-clusters, with only few $T_\lambda$ edges between them. For $\lambda_1<\lambda_2$, we get the tree $T_{\lambda_2}$ from $T_{\lambda_1}$ by contracting the edges with labels in $(\lambda_1,\lambda_2]$. Thus, if we have the collection of $\lambda$-clusters for all $\lambda\in(-\infty,\infty)$, then we can reconstruct the trees $T_\lambda$, and hence the macroscopic structure of the NCE seems to determine that of the MST. However, this is only an intuitive description: we ignored that there are a lot of small $\lambda$-clusters for any $\lambda$, hence it is not at all clear that the NCESL still determines an object that can be the scaling limit of the MST. We will see in the next section how one can build an actual proof.

\begin{wrapfigure}[5]{r}{1.1 in}
\begin{center}
\vspace{-0.3 in}
\epsfxsize=1 in \epsffile{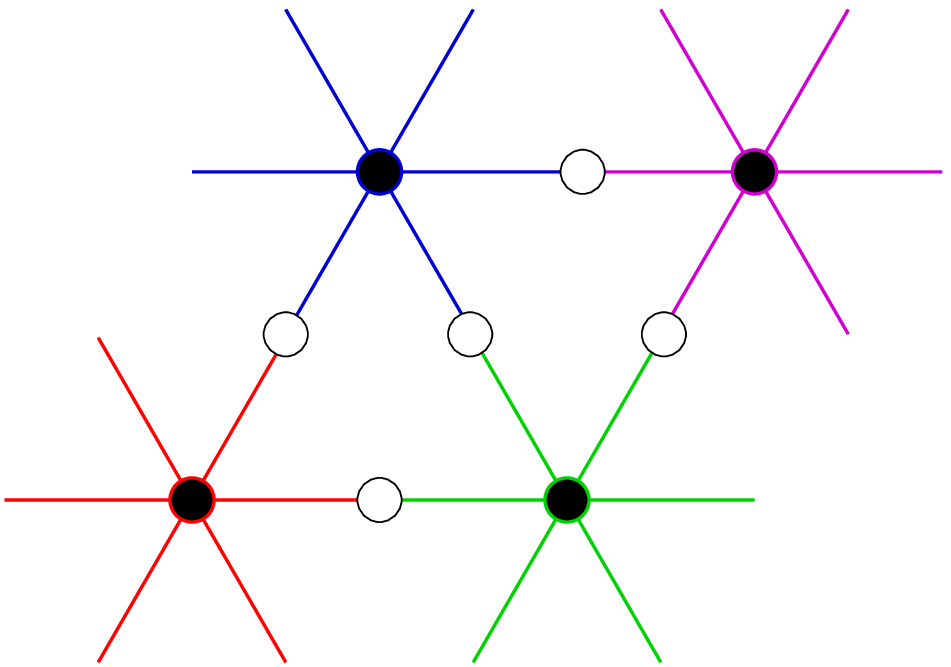}
\end{center}
\end{wrapfigure}

Since we have a proof of the existence and properties of the NCESL only for site percolation on the triangular lattice $\Delta$, if we want to use this to build the MST scaling limit, we will need a version of the MST that uses Unif$[0,1]$ vertex labels $\{V(x)\}$ on $\Delta$. So, replace each edge of $\Delta$ by two in series, and for each new edge $e$, denote its endpoint that was originally a vertex of $\Delta$ by $e^*$.  Then, let $U(e):=V(e^*)$. The MST using these edge labels $\{U(e)\}$ will inherit the right connectivity properties from the percolation ensemble. Our strongest results will apply to this model, but some of them will also hold for subsequential limits of the usual MST on $\Z^2$, known to exist by \cite{\ABNW}.

We note that the MST is also the union of the invasion trees of {\bf Invasion Percolation}, see 
\cite{\Alexander,\ABNW,\LPSMSF} and the references there.

\section{How do we see the MST in the scaling limit?}\label{s.MSTSL}


\begin{wrapfigure}[4]{r}{1.2 in}
\begin{center}
\vspace{-0.9 in}
\epsfysize=1.2 in \epsffile{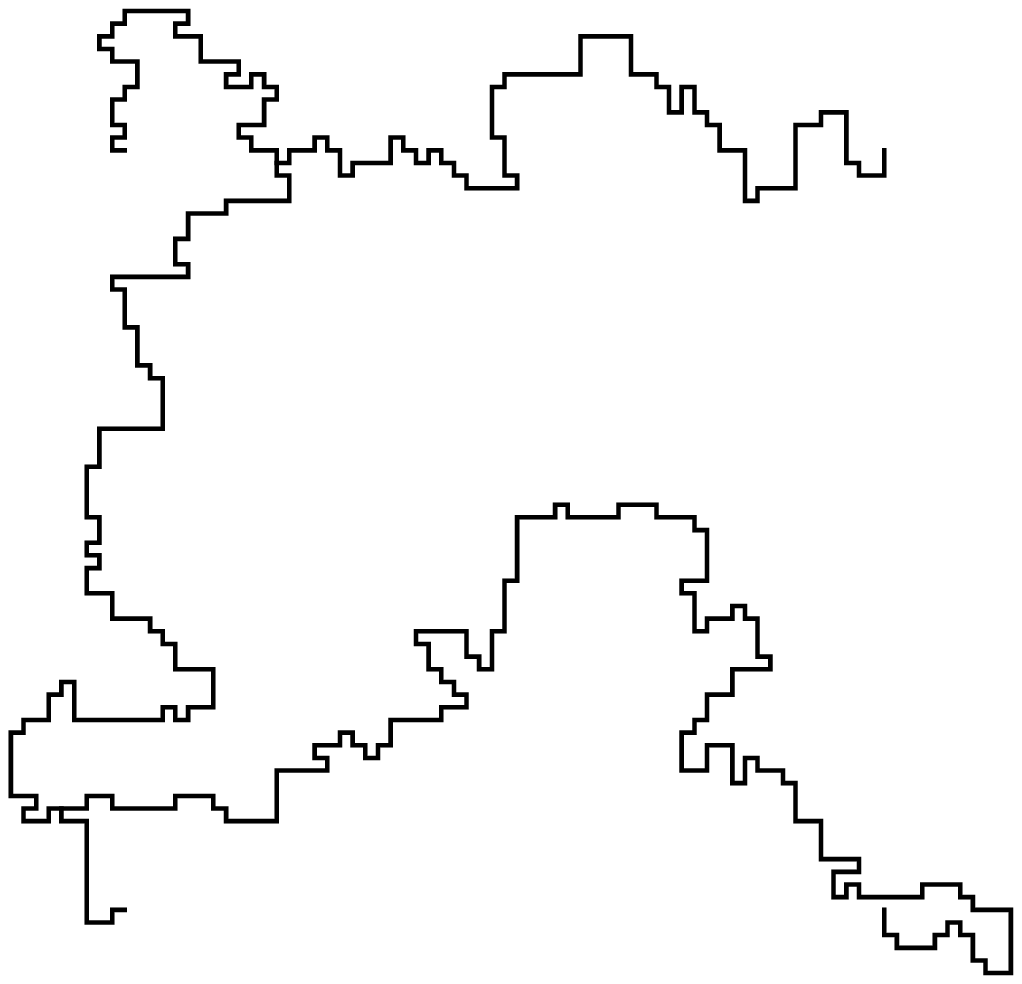}
\end{center}
\end{wrapfigure}

The MST scaling limit (MSTSL) should be a random spanning tree on all the points of the real plane. This, of course, cannot be a random subset of the plane, even though the discrete MST is a subgraph of the lattice. Rather, for any pair of points, $x,y$, we want the MSTSL path $P(x,y)$ connecting them. (And this collection of paths should satisfy some compatibility relations: the symmetric difference of $P(x,y)$ and $P(x,z)$ should be $P(y,z)$, and more generally, the subtree of the MSTSL connecting a finite set of given points should be part of the subtree connecting a larger finite set.) In fact, we fix a countable dense set of points $Z$ in the plane, and the MSTSL will be a collection of paths connecting all pairs of $Z$. This way, almost sure results for the path between a fixed pair will hold for all the pairs in $Z$ simultaneously. Then, for arbitrary points $x,y$, we can take sequences $x_n\to x$ and $y_n\to y$ with $x_n,y_n\in Z$, and take the limit of the paths $P(x_n,y_n)$ (with the metric given by the infimum of $L^\infty$-distances over all possible parameterizations of the paths). There will be pairs $x,y$ for which there are at least two different limiting paths; in this sense, the MSTSL is not exactly a tree. Nevertheless, this seems to be the best possible notion of the MST scaling limit. See \cite{\ABNW} for the precise description and basic topological properties of the subsequential scaling limits of the MST.

We now show that the MST path joining $x,y\in\R^2$ is determined by the NCE in such a way that the MSTSL will be a measurable function of the NCESL.

Fix some $\lambda_1 \ll 0$, so that even the outermost $\lambda_1$-clusters are small. As explained in Section~\ref{s.MST}, it is enough to find the path between the outermost clusters of $x,y$ in the cluster tree $T_{\lambda_1}$ in the scaling limit, then $\lambda_1\downarrow -\infty$ will give the MSTSL. 

\begin{wrapfigure}[8]{r}{1.8 in}
\begin{center}
\vspace{-0.36 in}
\SetLabels
\R(.1*.58)$x$\\
\R(.66*.52)$y$\\
\endSetLabels
\centerline{
\AffixLabels{%
\psfig{height=1.4in,figure=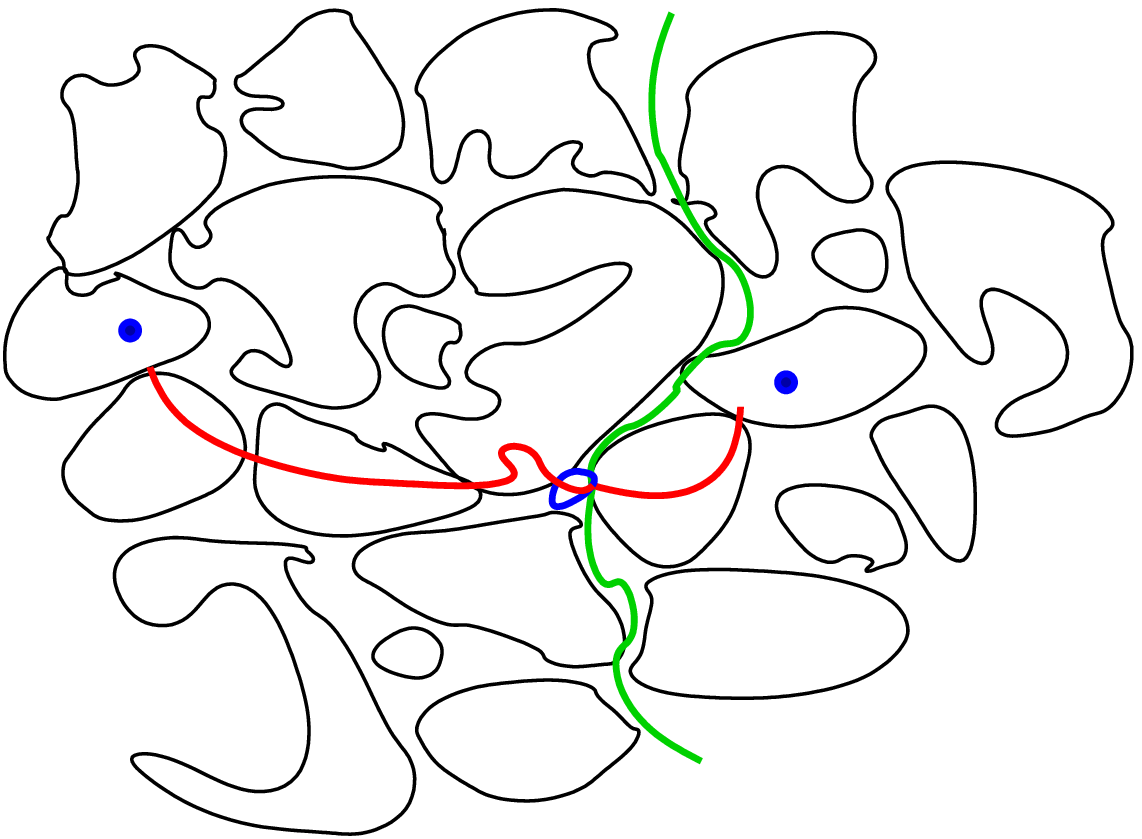}
}}
\end{center}
\end{wrapfigure}

The basic difficulty now is that there are infinitely many outermost $\lambda_1$-clusters even in a bounded region, so the MST on this cluster-graph may not be well-defined. So, take only the outermost $\lambda_1$-clusters with {diameter $\geq \eps$}.

Using RSW arguments, one can show that if $\eps$ is small enough, then there is a finite path between $x$ and $y$ in this cluster-graph. The resulting MST path uses labels $\leq \lambda_2$ only, since we need to open only finitely many cut-edges.

Now take $\delta \ll \eps$, and the corresponding new cluster graph. We claim that if $\eps$ was small enough, then 
the MST path in the new $\delta$-cutoff cluster graph is the same. 

Since the old path is still available in the $\delta$-cutoff cluster graph, the new path also has labels $\leq \lambda_2$. 
If the path goes through an outermost $\lambda_1$-cluster that has diameter $\delta$, then there is a $\lambda_1$-cut-edge $e$ on the path that is only $\delta$-important. But then, moving all the cut-edge labels from at most $\lambda_2$ to $\lambda_1$, this $e$ becomes very important at level $\lambda_1$, although it was very little important before. This contradicts the stability result (2) of Section~\ref{s.DPNCE}. Therefore, the  $T_{\lambda_1}$ path does not go through very small $\lambda_1$-clusters, so it is visible also in the scaling limit, and we are done.

\begin{wrapfigure}[4]{r}{2 in}
\begin{center}
\vspace{-0.45 in}
\SetLabels
\R(.14*.42)$x$\\
\R(.66*.55)$y$\\
\R(.48*.19)$e$\\
\endSetLabels
\centerline{
\AffixLabels{%
\psfig{height=1in,figure=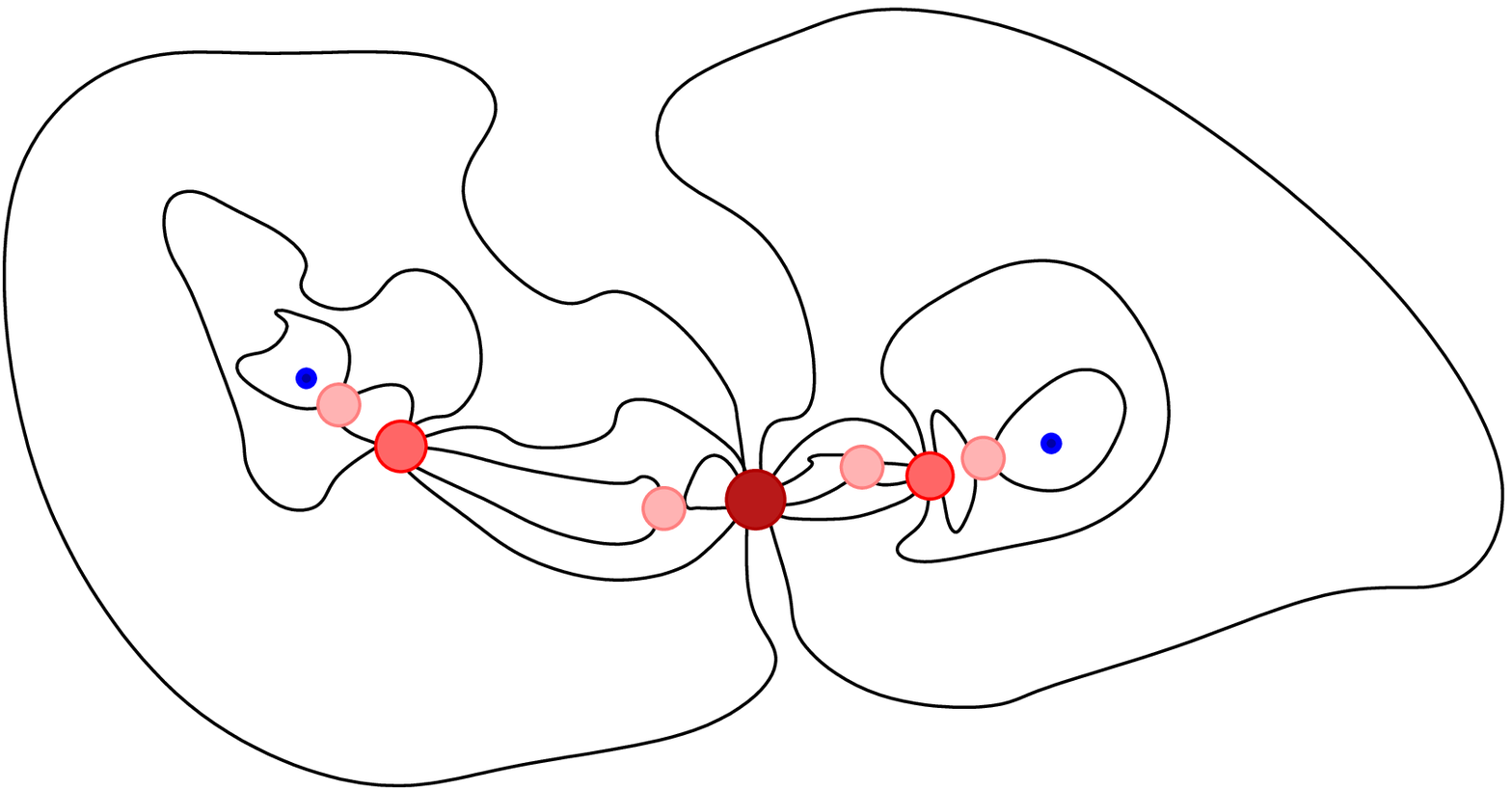}
}}
\end{center}
\end{wrapfigure}

Alternatively, one may try to find the minimal level where $x$ and $y$ become connected, at the cut-edge $e_{x,y}$, then repeat this between $e_{x,y}$ and $x$ and $y$ dyadically, until the entire path is recovered. However, when the mesh $\eta\downarrow 0$, the labels on the path near $x$ and $y$ blow up, hence $e_{x,y}$ in the scaling limit ``simultaneously coincides'' with both $x$ and $y$, so the procedure does not make sense. Hence we again need some macroscopic cut-off $\eps\downarrow 0$, for which proving the  convergence seems harder than above.

\section{Topology of the MST scaling limit}\label{s.topology}

In \cite{\ABNW} it was proved for any subsequential limit that almost all vertices are leafs and that there is a uniform bound on the vertex degrees. We can now prove stronger results: {\it In either lattice,  there are no degree $\geq 5$ points. For any two points in the plane, the MSTSL path joining them is a.s.~unique and simple (not even a figure of 6).} However, we do not know, for instance, if there are degree 4 points.

\begin{wrapfigure}[11]{r}{1.3in}
\begin{center}
\vspace{-0.35 in}
\epsfxsize=1.3 in \epsffile{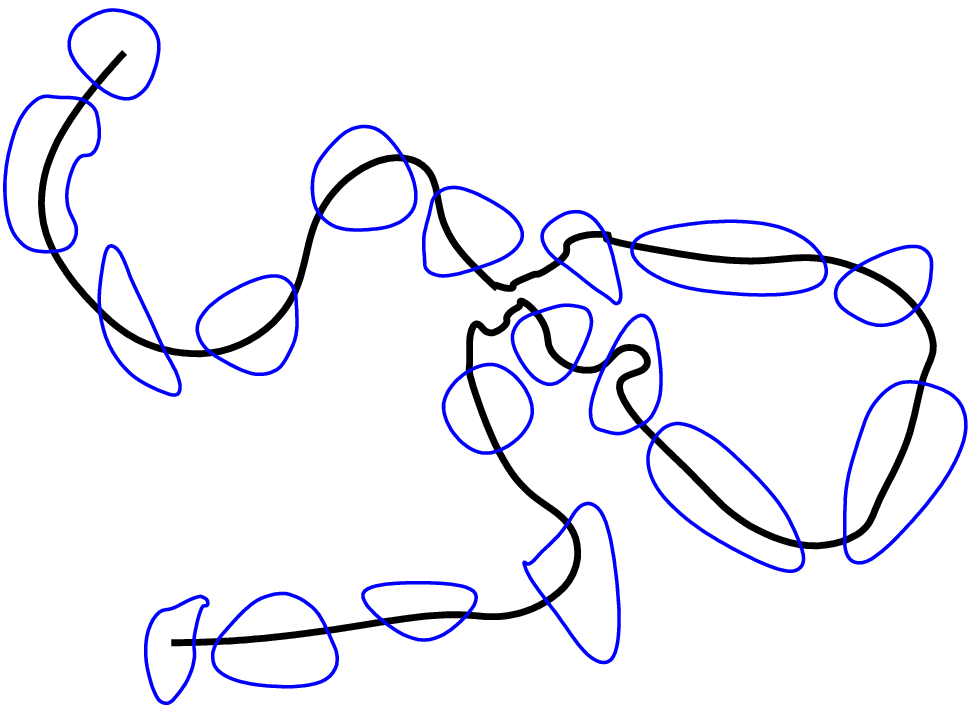}
\vskip 0.1 in
\SetLabels
\R(.19*.92)$x$\\
\R(.13*.01)$y$\\
\R(.54*.65)$z$\\
\endSetLabels
\centerline{
\AffixLabels{%
\epsfxsize=1.2 in \epsffile{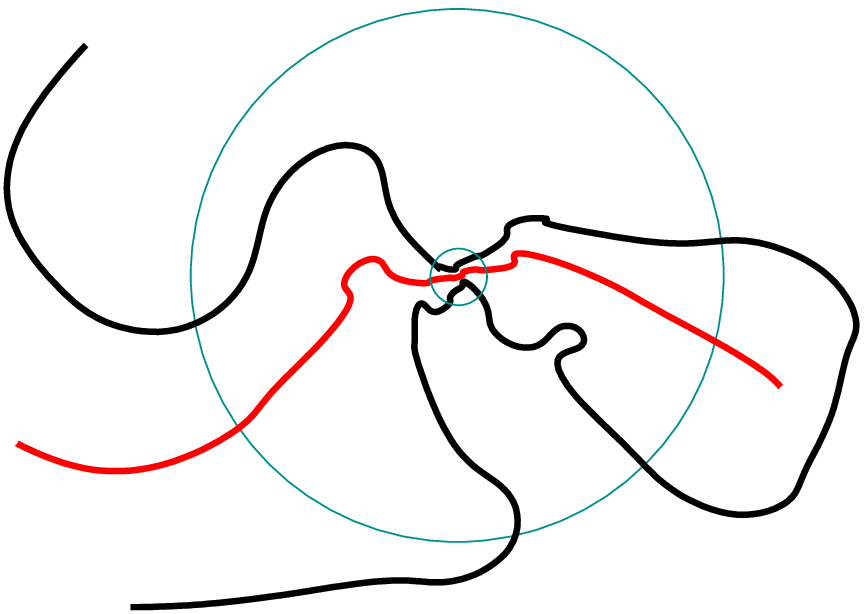}
}}
\end{center}
\end{wrapfigure}

Here is a proof sketch of the simple path claim. Consider a generic nearly non-simple path between $x$ and $y$. Take $\lambda_1$ very negative, and consider the part of $P(x,y)$ between the outermost $\lambda_1$-clusters surrounding $x$ and $y$. This part will be entirely below some finite level $\lambda_2$. On the other hand, since the $\lambda_1$-clusters all have small diameters, on the path there are labels above $\lambda_1$ ``all over the place''. But this implies that there must exist two macroscopic dual arms with labels all above $\lambda_1$ that force $P(x,y)$ go around the almost touch-point $z$. Altogether, we have 4 primal arms below level $\lambda_2$, and 2 dual arms with labels above $\lambda_1$, i.e., a six arm event around $z$ within a $W_\lambda$-modification, as in Claim (2) of Section~\ref{s.DPNCE}. Since it is known even on $\Z^2$ that the six-arm event does not happen at criticality, touch-points are ruled out.

\section{Conformal non-invariance?}\label{s.noninv}

\begin{wrapfigure}[6]{r}{1.3in}
\begin{center}
\vspace{-0.65 in}
\epsfysize=1.3 in \epsffile{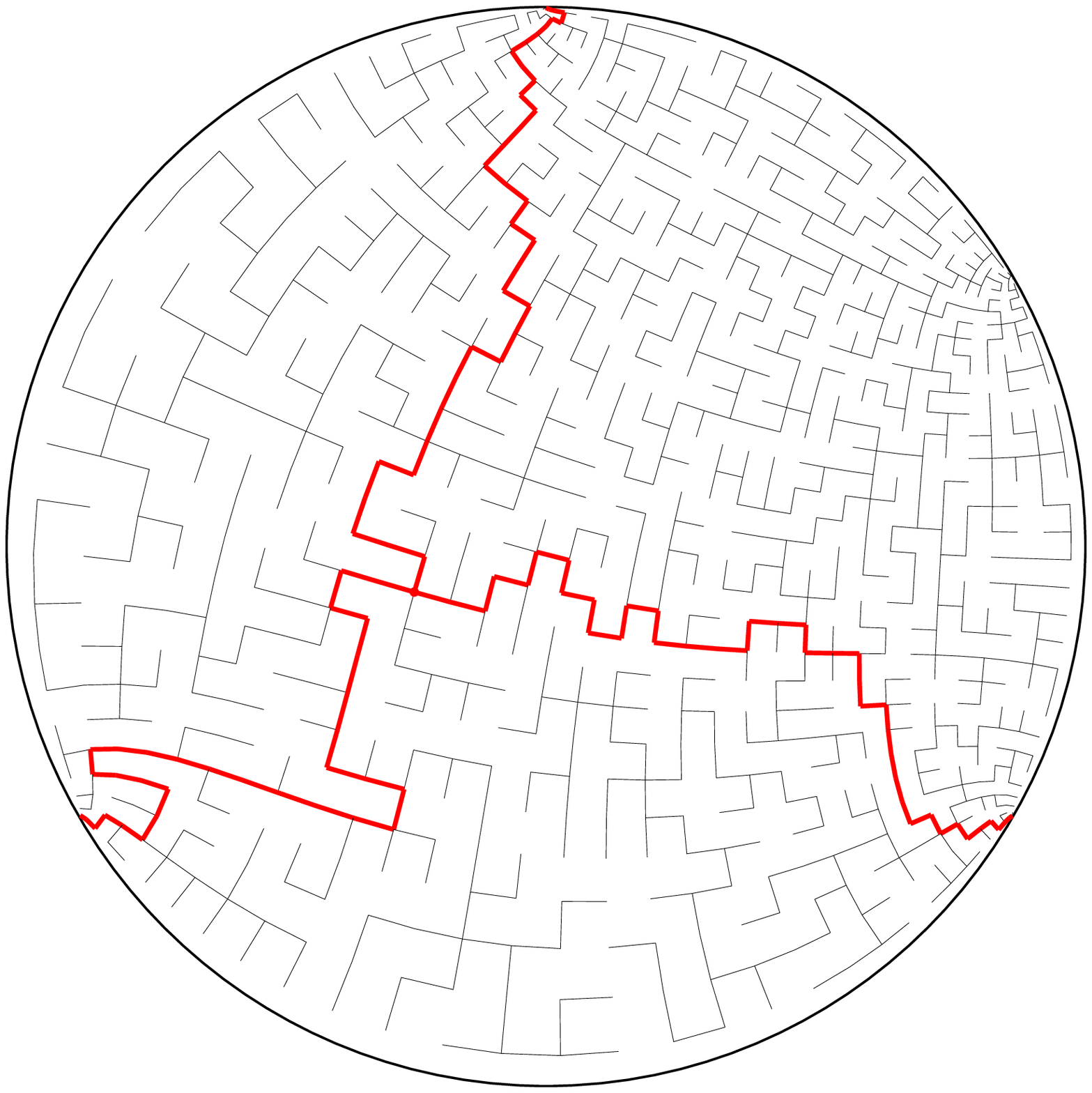}
\end{center}
\end{wrapfigure}

Scaling and rotational invariance of MSTSL follows easily from the same properties of NCESL, but the conformal {\it co}variance of the latter suggests {\bf conformal non-invariance of the MSTSL}. This is also supported by careful simulations \cite{\RGB}: the law of the trifurcation point on the first picture (a conformal image of a discrete square) is not invariant under rotation by $2\pi/3$. Here is a simplified version of what should be proved, leaving as an exercise to figure out the exact connection to the problem: 

%


Take an $n\times n$ square, with Unif$([0,1/5] \cup [4/5,1])$ labels on the left half, and 
Unif$[2/5, 3/5]$ labels on right. Take the MST path between the endpoints $x,y$ of the vertical half-line, 
and consider its segment between the $\eps n$-neighbourhoods of $x$ and $y$ (with high probability, there is only one). 
Does this segment feel the asymmetry between the two halves, as $n\to\infty$? It would be enough to prove, e.g., that the probability that it is contained in this or that side has different limits as $n\to\infty$. (In fact, my guess is that these two limit probabilities converge to 0 as $\eps\downarrow 0$ with different exponents in $\eps$.) It would also suffice to show that the limit law of this path is different from the law in a symmetric $n\times n$ square.

\begin{wrapfigure}[6]{r}{1.3in}
\begin{center}
\vspace{-0.52 in}
\epsfysize=1.2 in \epsffile{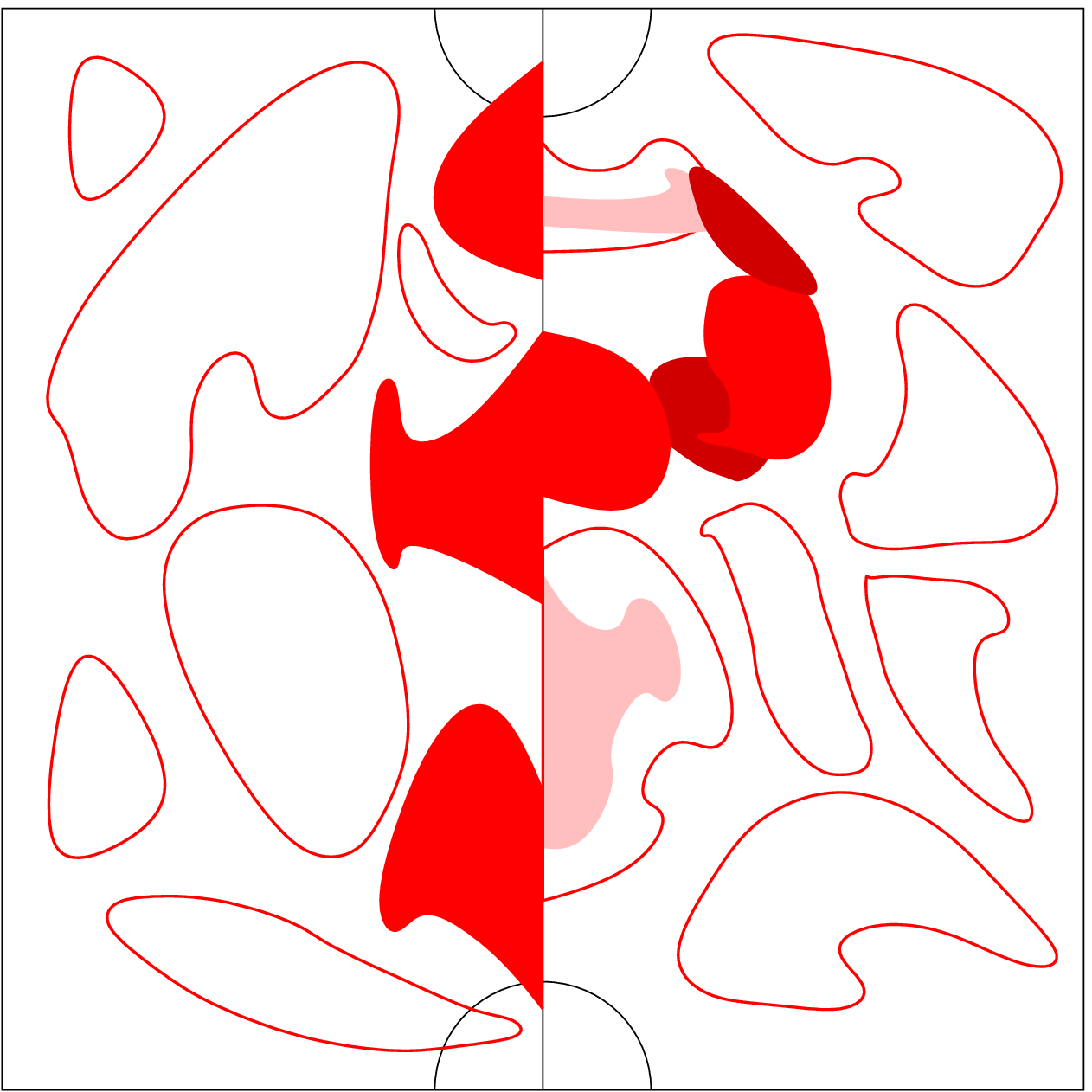}
\end{center}
\end{wrapfigure}

Now, what is the ``obvious asymmetry''? The edges open at level 1/2 form critical percolation, hence, for small $\eps$, it is unlikely that there is a 1/2-cluster connecting the two neighbourhoods. But then, the invasion tree from $x$ to $y$ will certainly use edges on the right side to travel between the 1/2-clusters, since all the labels larger than 1/2 are smaller on the right than on the left. On the other hand, the invasion tree will explore the entire 1/2-cluster on the left once having entered one, since all the labels smaller than 1/2 are smaller on the left than on the right, while it will explore only parts of the 1/2-clusters on the right. So, the MST path has more options on the right side for long distances, while more options on the left for short distances. The effects of this competition have to be understood well in order to produce an actual proof, but simulations suggest that the path spends more time on the right than on the left, and seems to intersect the midline less than in the symmetric situation.

This conformal non-invariance proof might be easier for Invasion Percolation, using known differences between invasion and critical percolation clusters  \cite{\DamronSV}.

\bibliographystyle{ws-procs975x65}
\bibliography{ws-pro-sample}

\end{document}